\documentclass[a4paper,10pt]{article}
\usepackage[utf8]{inputenc}
\usepackage[english]{babel}

\usepackage{amssymb}\usepackage{amsmath}
\usepackage{graphics}
\usepackage{epsfig}
\usepackage{bm}
\def\bl{\rule[-1mm]{2.4mm}{2.4mm}}

\def\be{\begin{equation}}
\def\ee{\end{equation}}

\newtheorem{thrm}{\bf Theorem}
\newtheorem{lmm}{\bf Lemma}
\newtheorem{dfn}{\bf Definition}
\newtheorem{rmk}{\bf Remark}

\begin{document}

\title {Conformal mapping of rectangular heptagons II}
\author{A.B.~Bogatyrev
\and  O.A.~Grigor'ev\thanks{Supported by RSCF grant 16-11-10349}}
\date{}
\maketitle
\abstract{
A new analytical method for the conformal  mapping of a rectangular polygons with 
a straight angle at infinity to a half plane and back is proposed. The method is 
based on the observation that SC integral in this case is an abelian integral on a hyperelliptic curve,
so it may be represented in terms of Riemann theta functions. The approach is illustrated 
by the computation of 2D flow of ideal fluid above rectangular underlying surface
and the computation of the capacities of multi component rectangular condensers with axial symmetry.}



\section{Introduction}
Numerical methods designed to evaluate a conformal mapping from a certain "standard domain" (say, an upper half-plane) are as 
numerous as are their applications in the area of computational physics \cite{TD}. For some domains only a finite set of parameters 
have to be found to substitute it into a known analytical expression for a conformal map. Delivering an 
example of such type of domains, rectangular polygons are the main subject of this work.

Conformal mapping from the upper half-plane to an arbitrary simply connected polygon is given by the 
Christoffel-Schwarz (CS) integral. If a polygon is rectangular (i.e., its angles are integer multiples of $\frac{\pi}{2}$), 
the CS integral is an abelian integral on a certain hyperelliptic surface. This observation  
along with the well-developed function theory on Riemann surfaces, was crucial for a novel method \cite{AB,BHY,Gr} 
of evaluation of conformal maps from an upper half-plane to a rectangular polygon and back. 

In a rectangular polygon we encounter four possible junctions of the neighboring edges. They may intersect  (1) at a right angle
(possibly intruding), (2) at a full angle $2\pi$ (a cut) (3) at zero angle (seminfinite "channel") and (4) at a straight angle $\pi$ at infinity. The first three cases were considered in detail in the papers \cite{AB, Gr}, however the remaining case: angle $\pi$ at infinity 
is also practically important and requires some modification of the technique used. 

As an evolution of the paper \cite{AB}, we consider polygons containing six finite vertices and an additional one at infinity. 
The conformal map in question as well as its inverse can be expressed by means of theta functions on a curve. For the latter 
there exist a robust and effective method of computation \cite{DHB} with controllable accuracy. Therefore we can guarantee the 
machine accuracy for  the conformal mapping uniformly in the polygon/halfplane.

In Sect. 2 we specify the spaces of polygons we work with, next section reminds several facts about Riemann surfaces 
and theta functions relevant to our study. Section 4 describes two moduli spaces of curves we deal with in what follows. 
The computational algorithm is formulated in Sect. 5. In the concluding Sect. 6 we apply the obtained algorithms to the 2D flow of ideal fluid above the rectangular underlying surface ("city landscape``) and to the  evaluation of the capacity of a planar rectangular condensers with axial symmetry.

\section{Spaces of polygons}
\begin{figure}[h]
\begin{picture}(170,65)
\thicklines
\put(5,15){\line(1,0){20}}
\put(25,15){\line(0,1){20}}
\put(25,35){\line(1,0){15}}
\put(40,35){\line(0,-1){30}}
\put(40,5){\line(1,0){10}}
\put(50,5){\line(0,1){50}}
\put(50,55){\line(1,0){20}}
\put(26,12){$w_1$}
\put(26,22){$H_1$}
\put(26,32){$w_2$}
\put(31,37){$H_2$}
\put(35,32){$w_3$}
\put(34,18){$H_3$}
\put(35,3){$w_4$}
\put(43,1){$H_4$}
\put(51,3){$w_5$}
\put(51,28){$H_5$}
\put(51,52){$w_6$}
\put(10,72){$a)$}
\put(75,45){\line(1,0){50}}
\put(125,45){\line(0,-1){20}}
\put(125,25){\line(-1,0){25}}
\put(100,25){\line(0,1){35}}
\put(100,60){\line(-1,0){15}}
\put(85,60){\line(0,-1){50}}
\put(85,10){\line(1,0){60}}
\put(80,72){$b)$}
\put(126,46){$w_1$}
\put(126,22){$w_2$}
\put(97,22){$w_3$}
\put(101,61){$w_4$}
\put(81,61){$w_5$}
\put(81,7){$w_6$}
\put(20,52){$Polygon$}
\put(110,52){$Polygon$}
\end{picture}

\caption{\small $a)$ Polygon from the space ${\cal P}_{236}$ $~~~~~~~~~~~~~~~b)$ Overlapping polygon from the space ${\cal P}_{123}$ }
\label{Heptagons}
\end{figure}
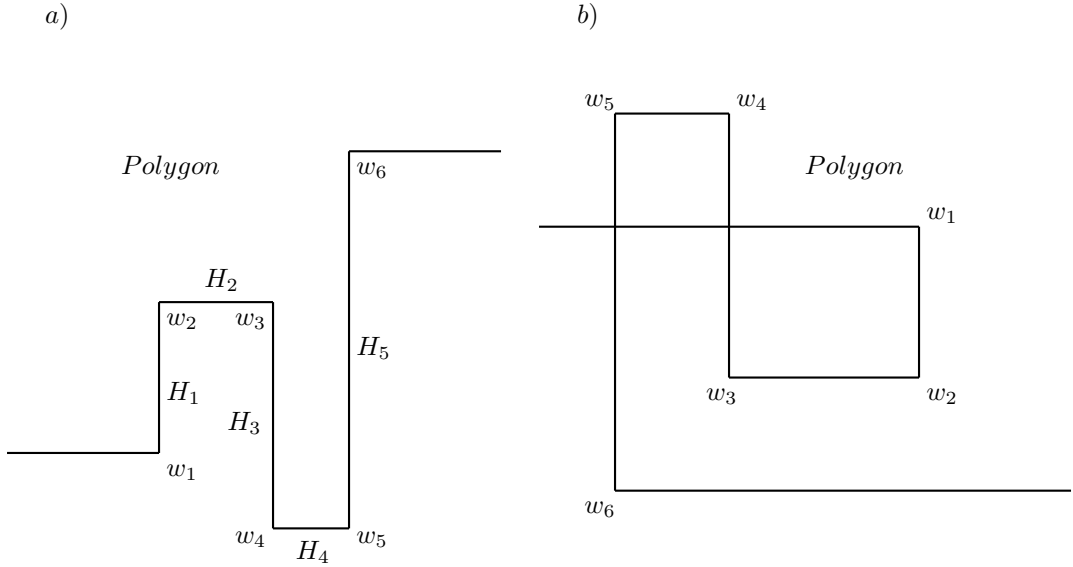

We consider simply connected heptagon with the straight angle at infinity, 
its sides are either vertical or horizontal. Exactly three of its right angles are intruding, that is equal to $\frac32\pi$ if measured inside
the domain. The combinatorial type of the polygon is given by the sequence of indexes  $0\le\sigma_1<\sigma_2<\sigma_3\le6$ of the intruding angles: 
the vertexes $\infty=:w_0,~w_1,\dots,~w_6$ of the polygon are naturally ordered  -- see Fig. \ref{Heptagons}. 
The space of heptagons  of a given combinatorial type $\sigma=(\sigma_1,\sigma_2,\sigma_3)$ we denote  
${\cal P}_\sigma$. We parametrize each space of polygons with the lengths of their sides, to which we ascribe signs for technical reasons:
\be
\label{HeptDim}
i^sH_s:=w_{s+1}-w_s,\qquad s=1,2,\dots,5.\\
\ee
The sign of the side length $H_s$ is negative iff $\sigma_1\le s<\sigma_2$ or $\sigma_3\le s$:
\be
\label{SignRule}
H_s\cdot P_\sigma(s+\frac12)<0,
\qquad P_\sigma(s):= \prod_{j=1}^3(s-\sigma_j).
\ee
Avoiding the degeneracy of polygons (such as vanishing isthmi) leads to further restrictions on the values $H_s$.
\be
\label{NoDegeneracy}
\begin{array}{c|cc}
\sigma= &{\rm Restriction}&\\
\hline
(126)&H_5-H_3<0& {\rm when}~H_4-H_2\le 0\\
(156)&H_1-H_3>0& {\rm when}~H_4-H_2\le 0\\
(123)&H_5<H_3~(<0)&\\
(456)&H_1>H_3~(>0)
\end{array}
\ee
Note that polygons of combinatorial types $\sigma=(123)$, $(234)$, $(345)$ and $(456)$
-- corresponding to three consecutive right turns of the boundary -- 
may be overlapping.

\begin{lmm}\label{H1H5}
The polygons of a given combinatorial type $\sigma$ make up a connected space ${\cal P}_\sigma$ of real 
dimension 5 with the global coordinates $H_1,H_2,H_3,H_4,H_5$ subjected to the sign rule (\ref{SignRule}) and inequalities
(\ref{NoDegeneracy}). ~~~\bl
\end{lmm}

A point from a space ${\cal P}_\sigma$ defines a polygon up to translations in the plane only. Those translations may be eliminated when necessary by the normalization e.g. $w_1:=0$. The reflection in the imaginary axis induces the mapping 
${\cal P}_{\sigma_1,\sigma_2,\sigma_3}\to$ ${\cal P}_{7-\sigma_3,7-\sigma_2,7-\sigma_1,}$ which in the intrinsic coordinates looks like 
$(H_1,H_2,\dots,H_5)\to$ $(-H_5,-H_4,\dots,-H_1)$.

\section{Hyperelliptic curves with six real branch points}
The conformal mapping of the upper half plane to any heptagon
from the space ${\cal P}_\sigma$ may be represented by the Christoffel-Schwarz integral. This integral can be lifted to a hyperelliptic curve
with six real branchpoints. In this section we briefly remind several facts about such curves. The survey is the short variant
of Sect. 2, 4 of \cite{AB} and does not contain any proofs, more details may be found in the textbooks like \cite{FK,GH,Nat}. 

\subsection{Algebraic model}
The double cover of the sphere with six real branch points
$x_1<x_2<\dots<x_5<x_6$
 is a compact genus two Riemann surface $X$ with the equation
 (of its affine part):
\be
\label{X}
y^2=\prod\limits_{s=1}^6(x-x_s),
\qquad  (x,y)\in\mathbb{C}^2.
\ee
 This curve admits a \emph{conformal involution} $J(x,y)=(x,-y)$ with six stationary points $p_s:=(x_s,0)$ and an \emph{anticonformal involution} (\emph{reflection})
$\bar{J}(x,y)=(\bar{x},\bar{y})$. The stationary points set of the latter has three components known as \emph{real ovals} of the curve. Each real oval is an embedded circle \cite{Nat} and doubly covers exactly one of the segments  $[x_2,x_3]$, $[x_4,x_5]$, $[x_6,x_1]\ni\infty$ of the extended real line $\hat{\mathbb R}:={\mathbb R}\cup\infty$. We denote those ovals as \emph{first, second and third} respectively. The  lift of the complimentary set of intervals
$[x_1,x_2]$, $[x_3,x_4]$, $[x_5,x_6]$ to the surface (\ref{X})
 gives us three \emph{coreal ovals} which make up the set of points fixed by another anticonformal involution ${\bar J}J=J{\bar J}$.

\subsection{Cycles, differentials, periods}\label{SectBasisCycles}
\begin{figure}
\includegraphics[scale=1.2]{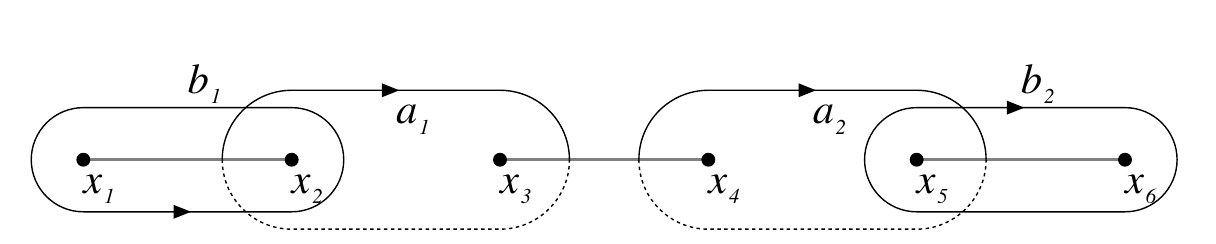}
\caption{\small Canonical basis in homologies of the curve $X$}
\label{BasisHomologies}
\end{figure}

We fix a special basis in the 1-homology space of the curve $X$ intrinsically related to the latter. The first and second real ovals give us two 1-cycles,  $a_1$ and $a_2$ respectively. Both cycles are oriented (up to simultaneous change of sign) as the boundary of a pair of pants obtained by removing real ovals from the surface. Two remaining cycles $b_1$ and $b_2$ are coreal ovals of the curve oriented so that the intersection matrix takes the canonical form --
see Fig. \ref{BasisHomologies}.

The reflection of the surface acts on the introduced basis as follows
\be
\label{abreflect}
\begin{array}{c}
\bar{J}a_s=a_s,
\quad
\bar{J}b_s=-b_s,
\end{array}
\qquad s=1,2.
\ee
Holomorphic differentials on the curve $X$ take the form
\be
\label{diffRep}
du_*=(C_{1*}x+C_{2*})y^{-1}dx,
\ee
with constant values $C_{1*}$, $C_{2*}$.
 The basis of differentials dual to the basis of cycles
\be
\int_{a_s}du_j:=\delta_{sj};
\qquad s,j =1,2,
\ee
determines Riemann period matrix $\Pi$ with the elements
\be
\label{periods}
\Pi_{sj}:=\int_{b_s}du_j;
\qquad s,j =1,2.
\ee
It is a classical result that $\Pi$ is symmetric and has positive definite imaginary part \cite{FK}.

From the symmetry properties (\ref{abreflect}) of the chosen basic cycles it readily follows that:
\begin{itemize}
 \item Normalized differentials are real ones, i.e. $\bar{J}du_s=\overline{du_s}$, in other words the coefficients $C_*$ in the representation (\ref{diffRep}) are real.
  \item Period matrix is purely imaginary,  therefore we can introduce the symmetric and positive definite real  matrix  $\Omega:=Im(\Pi)$,
\end{itemize}

\subsection{ Jacobian and Abel-Jacobi map }
\begin{dfn}
Given a Riemann period matrix $\Pi$, we define the full rank lattice
$$
L(\Pi)=\Pi\mathbb{Z}^2+ \mathbb{Z}^2 = \int_{H_1(X,\mathbb{Z})}du,
\qquad du:=(du_1,du_2)^t,
$$
in $\mathbb{C}^2$ and the 4-torus $Jac(X):=\mathbb{C}^2/L(\Pi)$ known as a Jacobian of the curve $X$.
\end{dfn}
This definition depends on the choice of the basis in the lattice $H_1(X,\mathbb{Z})$, other choices bring us to isomorphic tori.

It is convenient to represent the points  $u\in\mathbb{C}^2$ as a theta characteristic $[\epsilon,\epsilon']$, i.e. a couple of real 2-vectors (columns) $\epsilon, \epsilon'$:
\be
\label{ee'}
u=\frac12(\Pi\epsilon+\epsilon').
\ee
The points of Jacobian  $Jac(X)$ in this notation correspond to
two vectors with real entries modulo 2. Second order points of Jacobian are
represented as $2\times 2$ matrices with ${\mathbb Z}_2$ entries.
Our notation of theta characteristic as two column vectors is not commonly established: sometimes they use the transposed matrix.

\begin{dfn}
Abel-Jacobi (briefly: AJ) map is a correctly defined mapping from the surface $X$ to its Jacobian.
\be
\label{AJmap}
u(p):=\int_{p_1}^p du~~ mod~L(\Pi),
\qquad p_1:=(x_1,0); \quad du:=(du_1,du_2)^t,
\ee
\end{dfn}

From Riemann-Roch formula \cite{FK} it easily follows  that Abel-Jacobi map is a holomorphic embedding of the curve into its Jacobian. In Sect. \ref{theta} we give an explicit equation for the image of the genus two curve in its Jacobian. Let us meanwhile compute the images of the branching points $p_s=(x_s,0)$, $s=1,\dots,6$ of the curve $X$:

\centerline{
\begin{tabular}{c|c|c}
$p$ & $u(p)~ mod~ L(\Pi)$ & $[\epsilon, \epsilon'](u(p))$\\
\hline
$p_2$ & $\Pi^1/2$ & $\tiny
\left[\begin{array}{c} 10\\00\end{array}\right]$\\
$p_3$&$(\Pi^1+E^1)/2$&$\tiny
\left[\begin{array}{c} 11\\00\end{array}\right]$\\
$p_4$&$(\Pi^2+E^1)/2$&$\tiny
\left[\begin{array}{c} 01\\10\end{array}\right]$\\
$p_5$&$(\Pi^2+E^1+E^2)/2$&$\tiny
\left[\begin{array}{c}01\\11\end{array}\right]$\\
$p_6$&$(E^1+E^2)/2$&$\tiny
\left[\begin{array}{c}01\\01\end{array}\right]$\\
\label{AJPj}
\end{tabular}}

where $\Pi^s$ and $E^s$ are the $s$-th columns of the period and identity matrix respectively. One can notice that vector $\epsilon(u(p))$ is constant along the real ovals and $\epsilon'(u(p))$ is constant along the coreal ovals.

\subsection{Theta functions on genus two curves}\label{theta}
Here we give a short introduction to the theory of Riemann theta functions adapted to genus two surfaces. 
Three problems related to conformal mappings will be effectively solved in terms of Riemann theta functions:
\begin{itemize}
\item Localization of the curve inside its Jacobian;
\item Representation of the 2-sheeted projection of the curve to the sphere;
\item Evaluation of the normalized abelian integrals of the second and the third kinds which participate in the  decomposition of 
CS-integral.
\end{itemize}

\begin{dfn}
Let $u\in\mathbb{C}^2$ and $\Pi\in\mathbb{C}^{2\times2}$ be a Riemann
matrix, i.e. $\Pi=\Pi^t$ and $Im~\Pi>0$. The theta function of those two arguments is the following Fourier series

$$
\theta(u, \Pi):=\sum\limits_{m\in\mathbb{Z}^2}
\exp(2\pi im^tu+\pi i m^t\Pi m).
$$
Also they consider theta functions with characteristics which are the slight modification of the above theta.
$$
\theta[2\epsilon, 2\epsilon'](u, \Pi):=\sum\limits_{m\in\mathbb{Z}^2}
\exp(2\pi i(m+\epsilon)^t(u+\epsilon')+\pi i (m+\epsilon)^t\Pi (m+\epsilon))
$$
$$
=
\exp(i\pi\epsilon^t\Pi\epsilon+2i\pi\epsilon^t(u+\epsilon'))
\theta(u+\Pi\epsilon+\epsilon',\Pi),
\qquad \epsilon,\epsilon'\in\mathbb{R}^2.
$$
The matrix argument $\Pi$ of theta function is usually omitted when it is clear which matrix we mean.  Omitted vector argument $u$ is supposed to be zero and the appropriate function of $\Pi$ is called the theta constant:
$$
\theta[\epsilon, \epsilon']:=\theta[\epsilon, \epsilon'](0, \Pi).
$$
\end{dfn}
The convergence of these series grounds on the positive determinacy of
$Im~\Pi$.  The series has high convergence rate with well controlled accuracy \cite{DHB}.
Theta function has the following easily checked quasi-periodicity properties with respect to the lattice $L(\Pi):=\Pi\mathbb{Z}^2+\mathbb{Z}^2$:
\be
\label{quasiperiod}
\theta(u+\Pi m+m', \Pi)=\exp(-i\pi m^t\Pi m-2i\pi m^tu)\theta(u,\Pi),
\qquad m,m'\in\mathbb{Z}^2.
\ee
Quasi-periodicity of theta with characteristics is easily deduced from the last formula.

\begin{rmk}\label{RemTheta}
(i) Theta function with integer characteristics
$[2\epsilon, 2\epsilon']$ is either even or odd depending on the parity of the inner product $4\epsilon^t\cdot\epsilon'$. In particular, all odd theta constants are zeros.\\
(ii) It is convenient to represent integer theta characteristics as the sums of AJ images of the branch points, keeping only the indexes of those points: $[sk..l]$ means the sum modulo 2 of the theta characteristics of points
$u(p_s)$, $u(p_k)$, $\dots,u(p_l)$, e.g. $[35]$ corresponds to the characteristics
$\tiny \left[\begin{array}{c} 10\\11\end{array}\right]$.
\end{rmk}

\subsubsection{Image of Abel-Jacobi map}
The location of genus 2 curve embedded to its Jacobian may be reconstructed by solving a single equation.

\begin{thrm}[Riemann]
A point $u$ of Jacobian lies in the image $u(X)$ of Abel-Jacobi map
if and only if
\be
\label{JacIm}
\theta[35](u)=0.
\ee
\end{thrm}
For the image of AJ embedding localization for higher genus curves see \cite{AB3}.

\subsubsection{Projection to the sphere}\label{projection}
Any meromorphic function on the curve may be effectively calculated via the
Riemann theta functions once we know its divisor \cite{FK,Mum}. 
Take for instance the degree 2 function $x$ on the hyperelliptic curve (\ref{X}). This 
projection is unique if fully normalized. The incomplete normalization like $x(p_j)=0; x(p_0)=\infty$
brings us to a simple expression containing unspecified dilation
\be
\label{xofP}
x(p):= Const
\frac{\theta^2[jk35](u)}{\prod_\pm\theta[k35](u\pm u^0)},
\qquad u=u(p), 
\qquad u^0=u(p_0),
\ee
here $k$ is any index from the set $\{1,\dots,6\}$ different from $j$.

\subsubsection{Second and third kind abelian integrals}
On any Riemann surface there exists a unique abelian differential  of the third kind $d\eta_{pq}$ with two simple poles at 
the prescribed points $p,~q$ only, residues $+1,-1$ respectively and trivial periods along all $a-$ cycles. 
Differentiating with respect to the position of $p$, local coordinate $z$ at the pole being fixed, we get the $a-$ normalized  differential of the second kind $d\omega_{2p}=(z-p)^{-2}dz+holomorphic~form$  with a double pole at $p$. 
Note that the latter differential acquires a constant factor with the change of the local coordinate $z$ at the pole.  
Appropriate abelian integrals have closed expressions due to Riemann in terms of theta functions:

\be
\label{kind3}
\eta_{pq}(s):=\int_*^sd\eta_{pq}=
\log\frac{\theta[{\bm\epsilon}](u(s)-u(p))}
{\theta[{\bm\epsilon}](u(s)-u(q))},
\qquad p,q,s \in X;
\ee
\be
\omega_{2p}(s):=\int_*^sd\omega_{2p}=-
\nabla\log\theta[{\bm\epsilon}](u(s)-u(p))\cdot \frac{du}{dz}(p)
\label{kind2}
\ee
with any odd theta characteristic $[\bm\epsilon]$; the gradient is taken with respect to the vector argument of theta function.  
The tangent vector $du/dz(p)$ to the curve $X$ embedded to its Jacobian
is orthogonal to $\nabla\theta[35](u(p))$ according to \eqref{JacIm}, so the second kind integral in \eqref{kind2} may  be rewritten 
also as the determinant of a $2\times 2$ matrix.

\section{Two moduli spaces}
We are going to describe the space of CS integrals corresponding to the rectangular polygons from the spaces ${\cal P}_\sigma$. This space is an extension of the underlying space of genus 2 curves with real branch points.

\subsection{Genus 2 curves with three real ovals}
\begin{dfn} The moduli space  ${\cal M}_2\mathbb{R}_3$ is the factor  of   the cyclically ordered sextuples   of  points $(x_1, \dots, x_6)$ from
$\hat{\mathbb{R}}$ modulo the action of $PSL_2(\mathbb{R})$ (= real projective transformations of the equator $\hat{\mathbb{R}}$ 
conserving its orientation) .
\end{dfn}
Taking the algebraic model (\ref{X}) into account, the above definition of the moduli space is equivalent to the following \cite{AB}:
Moduli space ${\cal M}_2\mathbb{R}_3$ is the space of genus two Riemann surfaces with a reflection and 
three enumerated real ovals. Two surfaces are equivalent iff  there is a conformal $1-1$ mapping between 
them commuting with the reflections and respecting the marking of real ovals.

 Normalizing the branch points e.g. as $x_4=\infty$, $x_5=-1$, $x_6=0$, one gets  the global  coordinate system on the moduli space
 ${\cal M}_2\mathbb{R}_3$:
$$
0<x_1<x_2<x_3<\infty.
$$
Other normalizations bring us to other global coordinate systems in the same space. Yet another global coordinate system in this space is related to the periods of holomorphic differentials.

\begin{thrm}\cite{AB}
The period mapping $\Omega(X)$ from Sect. \ref{SectBasisCycles} is real analytic diffeomorphism from the moduli space ${\cal M}_2\mathbb{R}_3$  to the trihedral cone
\be
\label{PeriodsCone}
0<\Omega_{12}<min(\Omega_{11}, \Omega_{22}).
\ee
\end{thrm}

\begin{rmk}
The mappings from the period matrices to the (suitably normalized) set of branch points of the curve and back are effective. 
For genus two it is implemented by the Rosenhain formulas \cite{R} (see section \ref{theta}) in terms of theta constants.
\end{rmk}

\subsection{Curves with a marked point on a real oval}
The problems of conformal mapping use a slightly more sophisticated moduli space, that of the genus two real 
curves with a marked point on a real oval (see e.g.  \cite{AB2}).

\begin{dfn}
\label{M21R3}
By ${\cal M}_{2,1}\mathbb{R}_3$ we mean the sets of seven cyclically ordered  points $(x_0,x_1,\dots,x_6)$ in the real equator $\hat{\mathbb{R}}$ of the Riemann sphere modulo the action of $PSL_2({\mathbb R})$.
\end{dfn}
An argument similar to that in the previous subsection shows that we can introduce an equivalent but more geometric definition of 
the latter space as the space of genus two Riemann surfaces with three enumerated real ovals and a marked point $p_0\neq Jp_0$ on the third  oval. Two surfaces are equivalent iff there is a conformal mapping between them commuting with the reflections and respecting the enumeration of the ovals as well as the marked point. In the above definition \ref{M21R3}  the projection of the marked point to the sphere is designated as $x_0$, other coordinates $x_1,\dots,x_6$ are the projections of the branch points of the curve. The moduli space 
${\cal M}_{2,1}\mathbb{R}_3$ has natural projection to the space ${\cal M}_2\mathbb{R}_3$ consisting in forgetting of the marked point $x_0$.

One can introduce several coordinate systems on the space ${\cal M}_{2,1}\mathbb{R}_3$. Fixing three points of seven, say $x_4:=\infty$,
$x_5:=-1$, $x_6:=0$, the positions of the remaining four points of the 7-tuple  will give us the global coordinate system on the moduli space:
\be
0<x_0<x_1<x_2<x_3<\infty.
\label{moduli2}
\ee
Other normalizations of the 7-tuple of points result in different coordinate systems on the moduli space. It is a good exercise to show that the arising coordinate change is a real analytic 1-1 mapping from the 4D cone (\ref{moduli2}) to the appropriate cell.

Yet another coordinate system on ${\cal M}_{2,1}\mathbb{R}_3$ is the modification of that related to the periods matrix. Three variables
$\Omega_{11}$, $\Omega_{12}$, $\Omega_{22}$ are inherited from the space
${\cal M}_{2}\mathbb{R}_3$ and the fourth is either $u_1^0$ or $u_2^0$,
the component of the image of the marked point under AJ map.
The integration path for $u(p_0)$ is the interval of the third real oval
avoiding the branch point $p_6$.

\begin{lmm}\cite{AB}
Each mapping $(x_0,x_1,x_2,x_3) \to (\Omega_{11}, \Omega_{12}, \Omega_{22}, u_s^0)$, $s=1,2$ is a real analytic diffeomorphism of the cone (\ref{moduli2}) to the product of the cone (\ref{PeriodsCone}) and the interval $I/2$,   $I:=(0,1)$
\end{lmm}

\subsection{Christoffel-Schwarz differentials}\label{CSmap}
Consider the space ${\cal P}_\sigma$ as defined above. To each element of the moduli space
${\cal M}_{2,1}\mathbb{R}_3$ we ascribe the unique (up to dilation) abelian differential $dw_\sigma$, odd with respect to the 
hyperelliptic involution, double pole at the marked point $p_0$ and double zeros at three distinguished branchpoints $p_{\sigma_j}$,
$j=1,\dots,3$. In the algebraic model (\ref{X}) it takes the form:
\be
\label{dwab}
dw_\sigma=A\prod_{j=1}^3(x-x_{\sigma_j})\frac{dx}y,
\qquad A>0,
\ee
where $y^2=\prod_{j=1}^6(x-x_j)$ and the 7-tuple $\infty=:x_0,x_1,\dots,x_6$
represents an element of the space ${\cal M}_{2,1}\mathbb{R}_3$.

Each element of the space ${\cal P}_\sigma$ is the image of the upper half plane under the Christoffel-Schwarz map
\be
\label{CS}
w_{\sigma}(x):=\int_*^xdw_{\sigma}.
\ee
We have used a  natural embedding of the upper half plane $\mathbb{H}$ to the 
curve $X$ from the space ${\cal M}_2\mathbb{R}_3$ once there is a marked point on a real oval.
The set  $x^{-1}\mathbb{H}\subset  X$ has two components and exactly one of them has the marked point
$p_0$ on its boundary.

Slightly modifying the proof of a similar theorem considering rectangular heptagons with one zero angle vertex at infinity \cite{AB}, we obtain the following.

\begin{thrm}
\label{th2}
Christoffel-Schwarz mapping $w_{\sigma}(x)$
induces a real analytic diffeomorphic map from the decorated moduli space
${\cal M}_{2,1}\mathbb{R}_3\times\mathbb{R}^+$ to the heptagon space ${\cal P}_\sigma$.
\end{thrm}

Christoffel-Schwarz differential (\ref{dwab}) can be decomposed into a sum of elementary ones:
\be
dw_\sigma=c\cdot(d\omega_{2p_0}-J^*d\omega_{2p_0})+h\cdot d\eta_{p_0Jp_0}+c_1\cdot du_1+c_2\cdot du_2,
\label{diffdecomp}
\ee
where $d\omega_{2p}$ is an abelian differential of the second kind with a double pole at $p$, $d\eta_{pq}$ 
is an abelian differential of the third kind with poles at $p$ and $q$ only and residues $\pm1$; 
$du_1$ and $du_2$ are holomorphic differentials. All participating differentials are $a$-normalized and
$c$, $c_1,c_2$, $h$q are real constants. It remains to get the effective evaluation of the CS integral itself and all the auxiliary parameters.
Using Riemann's  formulas (\ref{kind3}), (\ref{kind2})  we represent the integral of (\ref{diffdecomp}) in terms of the jacobian variables
$$
w_\sigma(u)=c\det||\nabla_u\log\bigl(\theta[{\bm\epsilon}](u^0 - u)\theta[{\bm\epsilon}](u^0 + u)\bigl),~\nabla_u\theta[35](u^0)||
$$
\be
\label{wab}
+h\log\frac{\theta[{\bm\epsilon}](u^0 - u)}{\theta[{\bm\epsilon}](u^0 + u)}+ c_1u_1+c_2u_2.
\ee
where  $u:=(u_1,u_2)^t=u(p)$; $u^0:=u(p_0)$ and the periods matrix for theta is $i\Omega$.

\section{Algorithm of conformal mapping}
Based on the formulas of the previous section, we can propose the algorithm of conformal mapping of the heptagon to the half-plane and vice versa. First of all, given the heptagon we have to determine the corresponding point of the decorated moduli space ${\cal M}_{2,1}\mathbb{R}_3\times\mathbb{R}_+$.

\subsection{Auxiliary parameters}
Given coordinates of the polygon, we have to determine nine real parameters: the imaginary part $\Omega$ of the period matrix, 
the image $u^0:=u(p_0)$ of the marked point in the Jacobian of the curve, and 4 real scalars $c$, $c_1$ ,$c_2$, $h$. 
Those parameters give a solution to a system of nine real equations:

\be
\label{dweq0}
d\theta[35](u,i\Omega)\wedge
dw_\sigma=0, 
\qquad {\rm when}~u=u(p_{\sigma_j}), \quad j=1,2,3;
\ee
which means that CS differential $dw_\sigma$ has (double) zeros at three points $p_{\sigma_j}$;
\be
\label{ovalu0}
\theta[35](u^0,i\Omega)=0,
\ee
which means that the point $u^0$ lies on the AJ image of the curve in the Jacobian and finally
\be\label{Sides}
\begin{array}{l}
2H_1=4\pi c\theta_2[35](u^0,i\Omega)-4\pi h u_1^0 -c_1\Omega_{11}-c_2\Omega_{12};\\
2H_2=-c_1;\\
2H_4=c_2;\\
2H_5=4\pi c\theta_1[35](u^0, i\Omega)+4\pi h u_2^0 +c_1\Omega_{12}+c_2\Omega_{22};\\
H_1-H_3+H_5=\pi h
\end{array}
\ee
where $\theta_j[*](u,\dots):=\partial\theta[*](u,\dots)/\partial u_j$. Last five equations
specify the side lengths of the polygon and can be deduced by taking the CS integral along the real and coreal ovals 
and using the reciprocity laws \cite{FK,GH}. We note that they immediately give three unknown values $h$ and $c_1,c_2$ and 
linearly depend on the fourth $c$. The dependence of $\Omega$ and $u^0$ is not linear.

\begin{lmm}
Let the side lengths $H_1,H_2\dots,H_5$ satisfy the restrictions described in  Lemma \ref{H1H5} and the sign rule 
\eqref{SignRule}, then the system of nine equations \eqref{dweq0}, \eqref{ovalu0}, \eqref{Sides} has a unique solution 
with negative $c$, $2u^0\in I\times I$, $I=(0,1)$ and $\Omega$ in the trihedral cone \eqref{PeriodsCone}.
\end{lmm}
{\bf Proof sketch.} The existence and the uniqueness of the solution to (\ref{dweq0}) -- (\ref{Sides}) in the specified domain follows from the existence and the uniqueness of the conformal mapping of a given heptagon to the upper half plane. ~~~\bl

The strategy for the solution of the auxiliary set of equations \eqref{dweq0}-\eqref{Sides}  by parametric Newton method is discussed in  \cite{AB}.

\subsection{Mapping heptagon to the half plane and back}
Here we describe (somewhat sketchy, more details may be found in \cite{AB}) the algorithm of conformal mapping of a fixed heptagon 
from the space ${\cal P}_\sigma$ to the upper half- plane and back. First we solve the  system (\ref{dweq0}) -- (\ref{Sides})
of the auxiliary equations -- once for a given heptagon. Suppose a point $w^*$ lies in the heptagon normalized by the condition  $w_1=0$,
we consider a system of two complex equations
\be
\label{heptoH}
\begin{array}{r}
w_\sigma(u^*) = w^*,\\
\theta[35](u^*,i\Omega)=0,
\end{array}
\ee
with respect to the unknown complex 2-vector $u^*$. Following arguments in \cite{AB}, we assert that 
(\ref{heptoH}) has the unique solution $u^*$ with theta characteristic from the block
  $\tiny\left[
\begin{array}{cc}
-I~&I\\
-I~&I\\
\end{array}
\right]$, $\qquad I:=(0,1)$.
It is clear from the reflection principle for the conformal mappings that the set of two equations (\ref{heptoH}) may have many
solutions. We use theta characteristic to single out the unique one.

Substituting the solution $u^*$ to the right hand side of the expression (\ref{xofP}) we get the evaluation at the point $w^*$ of 
the conformal mapping $x(w)$ of the heptagon to the half plane with normalization $w_j$, $w_l$, $w_0\to$  $0,1,\infty$:
\be
\label{xofP2}
x^*=\pm\frac{\theta^2[lk35](u^0)}{\theta^2[lkj35]}
      \frac{\theta^2[jk35](u^*)}{\prod_\pm\theta[k35](u^*\pm u^0)},
\qquad k\neq j<l.
\ee
the sign $\pm$ in the latter formula is $(-1)^{\epsilon(l)\epsilon'(j)}$, where $[\epsilon(s),\epsilon'(s)]$ is the representation of the 
$u(p_s)$ as theta characteristic.
 
Conversely, given a point $x^*$ in the upper half plane we solve a system of two equations (\ref{xofP2}) and (\ref{JacIm}) with 
respect to a complex 2-vector $u^*$ with characteristics from the mentioned above block. 
Then substitute this solution to the formula (\ref{wab}) for CS integral in the jacobian variables
to get the image $w^*$ of the point $x^*$ in the heptagon.

\section{Applications}
One can enlarge the stock of heptagonal domains at the cost of minor complication of the involved moduli spaces and the algorithm.
Suppose three equations (\ref{dweq0}) are not satisfied and the zeros of the SC differentials have moved to the neighboring 
real or coreal ovals. Then the composite function  $w_\sigma(u(x))$ maps the upper half- plane to a heptagon with vertical or horizontal 
cuts emanating from the intruding right angles. The spaces of the decagons of this kind -- with six right angles, three full angles (spikes)
and one straight angle at infinity have dimension eight and the corresponding conformal mappings to the half- plane use the moduli space of real genus two curves with four marked points on their (co)real ovals. The parametric equations for the mapping itself are the same as above but they contain extra auxiliary parameters and equations for them. Now we have three more unknown points in the jacobian which gives six real values and six more real equations. Three of those mean the same as (\ref{ovalu0}): additional points live on the AJ images of (co)real ovals of the curve. Another three equations specify the lengths of three spikes. Note that six additional variables do not participate in the 
final parametric representation of the conformal mapping, however they influence the values of the involved auxiliary parameters.
We hope that an interested reader will easily reconstruct the mentioned set of 15 equations (much of those are linear) for 15 unknowns.

\subsection{A wind over the city} 
\begin{figure}[h]
\centering
\includegraphics[keepaspectratio=true,scale=0.4]{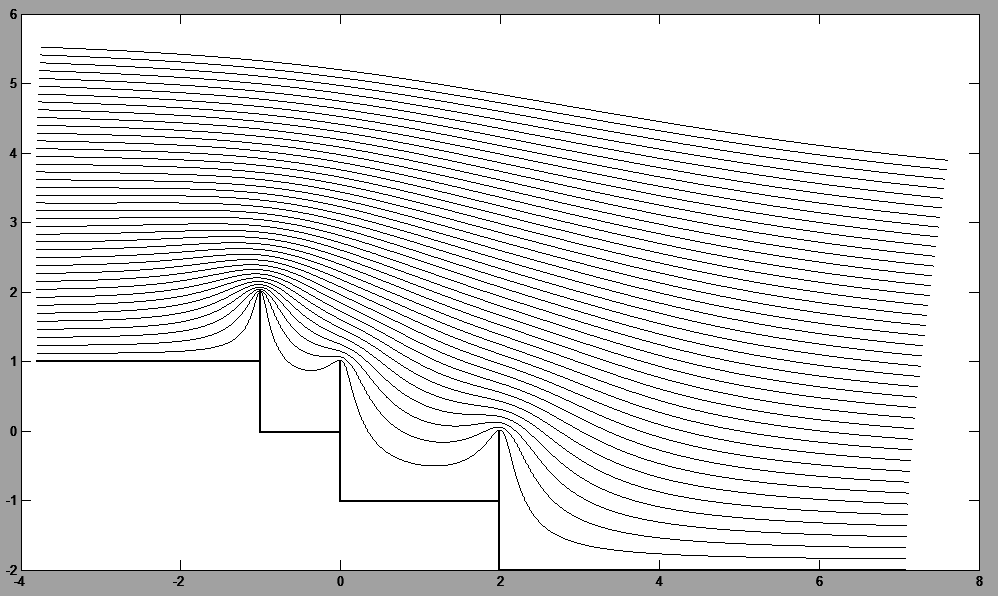}
\caption{Global behaviour of the stream lines.}
\label{channel1}
\end{figure}

\begin{figure}[h]
\centering
\includegraphics[keepaspectratio=true,scale=0.5]{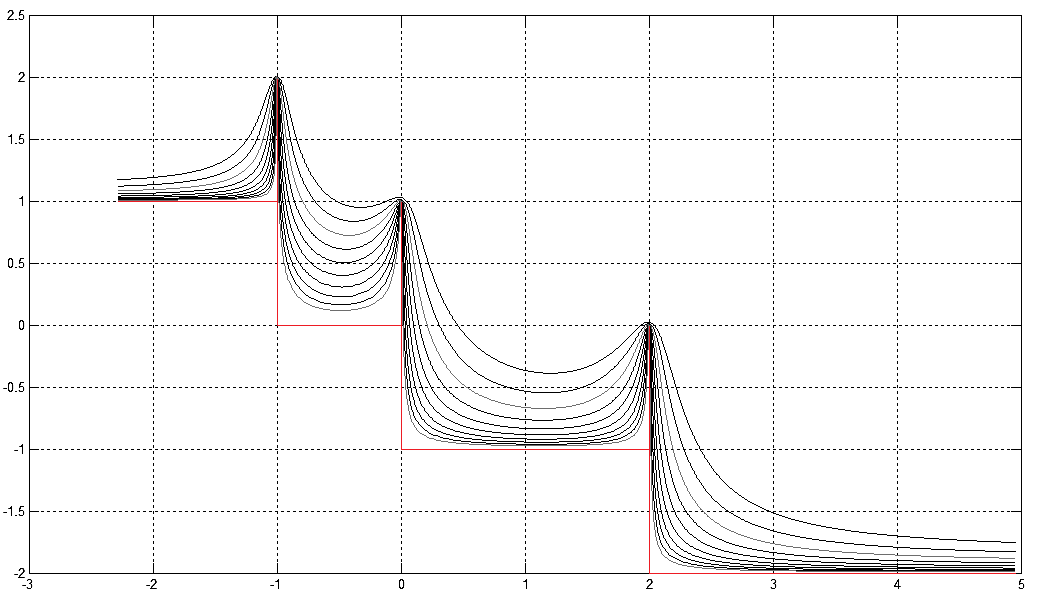}
\caption{Behaviour of the stream lines close to the border.}
\label{channel2}
\end{figure}

The stream lines of the ideal fluid flow above the flat surface are just the horizontal strait lines.
We conformally map the upper half plane to the decagon modeling the city landscape, with hydrodynamic normalization at infinity.
The horizontal lines are transplanted by this map to the streamlines of 2D ideal fluid above the surface with intricate geometry
as in Figs. \ref{channel1}, \ref{channel2} 

\subsection{Calculating capacities}
\begin{figure}[h!]
\centering
	\includegraphics[keepaspectratio=true,scale=0.5]{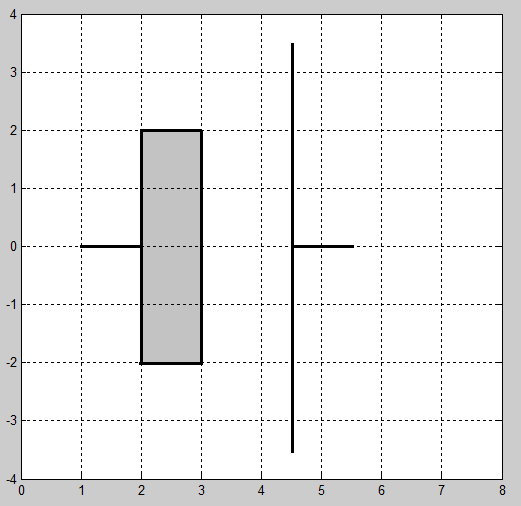}\hfill
	\includegraphics[keepaspectratio=true,scale=0.5]{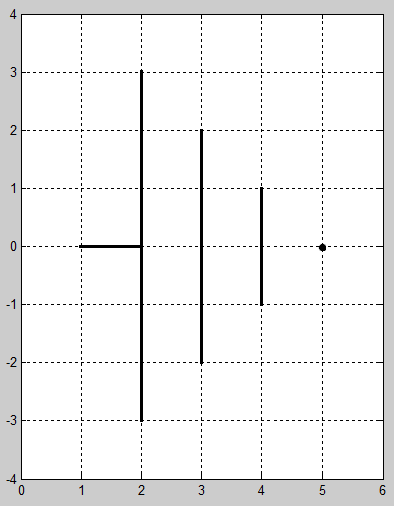}
	\caption{Left: Sign of a battery has capacity $C=2,3775929$;
	        Right: Capacity of the ground sign is $1,930955$.}
	\label{capacitors}
\end{figure}

Let $E$ being a bounded compact set (a condenser) in the complex plane. Its capacity $C:=Cap(E)$ is defined in terms of the asymptotic behaviour of Green's function $G_E(x)$  at infinity:
$$
G_E(x)=\log|x| -\log C +~o(1), \qquad x\to \infty.
$$
Once we may conformally map the exterior of one set to the exterior to the other with the hydrodynamical normalization 
$w(x)=x+O(|x|), \quad |x|\to\infty$, the capacities of two sets are equal. An exterior of an axisymmetric condenser
(with all the components on the axis) can always be mapped to the plane with several slots on a straight line. The capacity of the latter condenser may be found explicitly, say by Akhiezer or Sebbar-Falliero formulas. For the general formula of this kind condensers in terms of theta functions see \cite{BoGr}. In this way one can compute the capacities of the battery  and the ground symbols  
well known to every electrical engineer.  The dimensions of the condensers and their capacities are shown in the Fig. 
\ref{capacitors}.
 
\vspace{5mm}
\parbox{9cm}
{\it
119991 Russia, Moscow GSP-1, ul. Gubkina 8,\\
Institute for Numerical Mathematics,\\
Russian Academy of Sciences\\[3mm]
{\tt ab.bogatyrev@gmail.com, guelpho@mail.ru}}
\end{document}